\documentclass[12pt]{article}

\usepackage{amsmath}
\usepackage{lipsum}

\usepackage{amsthm}
\usepackage{mathrsfs}	
\usepackage{amssymb}
\usepackage{amscd}
\usepackage{setspace}	
\usepackage{tikz}		
\usepackage{tikz-cd}
\usetikzlibrary{matrix}	
\usepackage{graphicx}	
\graphicspath{ {pics/} }	
\usepackage{wrapfig}	
\usepackage{enumerate}	
\usepackage{bbm}		

\usepackage[utf8]{inputenc}	
\usepackage[english]{babel}	
\usepackage{blindtext}		
\usepackage[affil-it]{authblk}	

\usepackage{faktor}			
\usepackage{xparse}			

\usepackage{hyperref}		
\usepackage[lite]{amsrefs}

\newtheorem{theorem}{Theorem}[section] 

\newtheorem{definition}[theorem]{Definition}

\newtheorem{lemma}[theorem]{Lemma}

\newtheorem{corollary}[theorem]{Corollary}

\newtheorem{example}[theorem]{Example}

\newcommand{\acts}{\curvearrowright}
\newcommand{\asdim}{\text{asdim}}
\newcommand{\diam}{\text{diam}}

\newcommand{\Z}{\mathbb{Z}}
\newcommand{\C}{\mathbb{C}}
\newcommand{\N}{\mathbb{N}}
\newcommand{\dad}{\text{DAD}}

\begin{document}
\title{On Sharp Bounds for the Dynamic Asymptotic Dimension}
\author{Samantha Pilgrim}
\maketitle

\normalfont



\paragraph{Abstract:}We prove the dynamic asymptotic dimension of a free isometric action on a space of finite doubling dimension is either infinite or equal to the asymptotic dimension of the acting group; and give a full description of the dynamic asymptotic dimension of translation actions on compact Lie groups in terms of the amenability and asymptotic dimension of the acting group.  

\section{Introduction}
\normalfont
The dynamic asymptotic dimension (DAD) of a group action $\Gamma\acts X$ was first introduced by Guentner, Willett, and Yu \cite{dasdimGWY}, and was shown to relate to conditions used by Bartels, L\"{u}ck, and Reich in work on the Farrell Jones conjecture on manifold topology \cite{bartels} \cite{bartels_et_al}, as well as to the nuclear dimension of the crossed product \cite{winter2009nuclear} and calculations of its $K$-theory \cite{Guentner2016DynamicAD}.  The DAD is also related to other dynamical dimension theories \cite{kerr2017dimension} known to take into account both the asymptotic dimension of $\Gamma$ and the topological dimension of $X$.  However, it has been suspected since its introduction that the dynamic asymptotic dimension may actually coincide with the asymptotic dimension of $\Gamma$ whenever it is finite (it is possible to have $\dad(\Gamma\acts X) = \infty$ while $\asdim\Gamma<\infty$).  The case where $\Gamma$ is virtually cyclic is proved in \cite{amini2020dynamic}, providing an important test-case for this conjecture.  More generally, the work of \cite{warpedcones} gives an upper bound of $\asdim\Gamma + \dim X$ for isometric actions by finitely generated groups on manifolds; and \cite{conley2020borel} shows that many actions on Cantor sets have $\dad = \asdim \Gamma$.  

This note documents progress on a question \cite[8.6]{warpedcones} originally posed by Willett by giving the first sharp bounds for the dynamic asymptotic dimension of free isometric actions on spaces of dimension greater than (or equal to) zero, modulo the fairly mild assumption that $X$ has finite doubling dimension.  In this case we show the dimension is always either $\infty$ or $\asdim\Gamma$.  We accomplish this by introducing a new way of formulating the doubling dimension of a space, which we call Cantor decomposability.  This property, together with the theory of residually finite group actions, allows many isometric actions to be modeled by several sequences of partial dynamical systems on discrete spaces which asymptotically resemble the structure of $\Gamma$.  This property is reminiscent of box spaces of residually finite groups, and our investigations are partly inspired by past work on the asymptotic dimension of such spaces \cite{boxspacesDT}, as well as by the relation between box spaces and odometers \cite{elementaryamenableboxspaces}.  In line with \cite{boxspacesDT}, which describes among other things the asymptotic dimension of box spaces of residually finite groups; we use the theory of residually finite group actions introduced in \cite{kerr2011residually}, together with Cantor decomposability, to reduce the problem to geometry.  

We apply the main result in order to calculate the dimension of translation actions on compact Lie groups in terms of the amenability and asymptotic dimension of the acting group.  We also compute the dimension of many isometric actions by $\Z^n$.

\section{DAD of partial dynamical systems}

We begin by establishing some notation and terminology.  We will assume throughout that the group $\Gamma$ is finitely generated with finite generating set $F$, that $F = F^{-1}$ is symmetric, and that $e\in F$ (where $e$ denotes the identity in $\Gamma$).  The following definition comes from \cite[2.1]{Exel2015PartialDS}.  See \cite[Part I]{Exel2015PartialDS} for a more complete treatment of partial actions.  

\definition{A topological partial action of $\Gamma$ on a topological space $X$ is a pair $(\{D_\gamma\}_{\gamma\in \Gamma}, \{\theta_\gamma\}_{\gamma\in \Gamma})$ consisting of a collection $\{D_\gamma\}_{\gamma\in \Gamma}$ of subsets of $X$, and a collection $\{\theta_\gamma\}_{\gamma\in \Gamma}$ of homeomorphisms, $\theta_\gamma: D_{\gamma}\to D_{\gamma^{-1}}$ such that  

\begin{description}
\item(i) $D_e = X$, and $\theta_e$ is the identity map.  

\item(ii) $\theta_\gamma\circ \theta_\lambda\subseteq \theta_{\gamma\lambda}$, for all $\gamma$ and $\lambda$ in $\Gamma$.  
\end{description}

Here, the composition $\theta_\gamma\circ \theta_\lambda$ denotes the map whose domain is the set of all $x\in X$ for which $\theta_\gamma(\theta_\lambda(x))$ makes sense.  In other words, this is the set $\theta_\lambda^{-1}(D_\gamma) = \theta_{\lambda}^{-1}(D_{\lambda^{-1}}\cap D_\gamma)$.  The symbol ``$\subseteq$" means that the function on the right-hand-side is an extension of the function on the left-hand-side.  Notice also that $\theta_{\gamma^{-1}} = \theta_{\gamma}^{-1}$.  

A partial dynamical system is a quadruple $(X, \Gamma, \{D_\gamma\}_{\gamma\in F}, \{\theta_\gamma\}_{\gamma\in \Gamma})$.  An action $\Gamma\acts X$ is then a partial dynamical system with $D_\gamma = X$ for all $\gamma\in \Gamma$.  A partial dynamical system is free if $\theta_\gamma(x) = \theta_{\lambda}(x)$ iff $\gamma = \lambda$.  We may sometimes still write $\gamma\cdot x$ for $\theta_\gamma(x)$ when it is not ambiguous to do so.  

Some partial dynamical systems we will deal with later are restrictions of actions or partial actions, in which case we will write $\Gamma\acts A$ for a subset $A\subset X$ to denote the partial dynamical system $(A, \Gamma, \{A\cap h_\gamma^{-1}(A)\}_{\gamma\in F}, \{h_\gamma\}_{\gamma\in F})$ where $h_\gamma$ comes from $\Gamma\acts X$ or some partial action on $X$.  }

\definition{Let $(X, \Gamma, \{D_\gamma\}_{\gamma\in F}, \{\theta_\gamma\}_{\gamma\in \Gamma})$ and $(Y, \Gamma, \{E_\gamma\}_{\gamma\in F}, \{\rho_\gamma\}_{\gamma\in \Gamma})$ be topological partial dynamical systems.  By a conjugacy, we mean a homeomorphism $f: X\to Y$ such that $f(\theta_\gamma(x)) = \rho_\gamma (f(x))$ whenever $\theta_{\gamma}(x)$ is defined (so $\rho_\gamma(f(x))$ must be defined in this case).  }

\definition{Let $(X, \Gamma, \{D_\gamma\}_{\gamma\in \Gamma}, \{\theta_\gamma\}_{\gamma\in \Gamma})$ be a partial dynamical system and $S\subset \Gamma$ a finite subset with $S = S^{-1}$ and $e\in S$, and $A\subset X$.  An $S$-chain in $A$ is a finite sequence $x_0, \ldots, x_n$ of points in $A$ such that for all $0\leq i\leq n-1$, $\theta_{\gamma}(x_i) = x_{i+1}$ for some $\gamma\in S$.  Two points in $A$ are in the same $S$-component if they are connected by an $S$-chain in $A$.  We say a cover $\mathcal{V} = \{V_j\}_{j=0}^d$ is a $(d, S, M)$-cover for $\Gamma\acts X$ (or just $X$ if unambiguous) if all $S$-components of each $V_i$ have cardinality at most $M$.  }\label{dad notation}

\normalfont

\definition{The dynamic asymptotic dimension of a free partial dynamical system is the smallest integer $d$ such that for every finite subset $S\subset \Gamma$ with $S = S^{-1}$ and $e\in S$, there is $M>0$ and a $(d, S, M)$-cover for $(X, \Gamma, \{D_\gamma\}_{\gamma\in \Gamma}, \{\theta_\gamma\}_{\gamma\in \Gamma})$.  If no such $d$ exists, the dimension is defined to be $\infty$.  As we are assuming $\Gamma$ is finitely generated with finite generating set $F$, we can assume $S$ has the form $F^r$ for some $r>0$.  A $d$-dimensional control function for $(X, \Gamma, \{D_\gamma\}_{\gamma\in \Gamma}, \{\theta_\gamma\}_{\gamma\in \Gamma})$ (with $\Gamma$ finitely generated) is a function $D: \mathbb{N}\to \N$ such that for every integer $r>0$, there is a $(d, F^r, D(r))$-cover for $(X, \Gamma, \{D_\gamma\}_{\gamma\in \Gamma}, \{\theta_\gamma\}_{\gamma\in \Gamma})$.  }

\normalfont
We will use the notation described in the last two definitions (which relates to the definition of $\dad$ for free actions) even for actions which are not free.  

Notice that $\dad$ appears a priori to be sensitive to the topology of $X$ as the cover $\mathcal{U}$ is required to be open.  

\definition{If $(X_n, \Gamma, \{D^n_\gamma\}_{\gamma\in \Gamma}, \{\theta^n_\gamma\}_{\gamma\in \Gamma})$ is a sequence of partial dynamical systems, we say $\dad(\Gamma\acts X_n)_n\leq d$ uniformly if there is a single $d$ dimensional control function for all $n$.  }

\definition{Suppose $C_F(\Gamma)$ is a Cayley graph of a finitely generated group $\Gamma$ with $F$ being a generating set.  By this we mean the set $\Gamma$ with the metric given by $d(x, y) = |yx^{-1}|_F$ where $|w|_S$ denotes the minimal length of a word in $S$ which equals $w$ as an element of $\Gamma$.  Notice that this is invariant under right multiplication.  If $r>0$, an $r$-chain in $A\subset C_F(\Gamma)$ is a finite sequence $x_0, \ldots, x_n$ of points in $A$ such that $d(x_i, x_{i+1})\leq r$ for $0\leq i\leq n-1$.  Two points $g, h\in A\subset C_F(\Gamma)$ are in the same $r$-component of $A$ if they are connected by an $r$-chain in $A$.  We say a cover $\mathcal{U} = \{U_0, \ldots, U_d\}$ of $C_F(\Gamma)$ such that the $r$-components of each $U_i$ have diameter at most $M_r$ is a $(d, r, M_r)$-cover for $C_F(\Gamma)$.  The asymptotic dimension of $C_F(\Gamma)$, which we write as $\asdim C_F(\Gamma)$, is the least integer $d$ such that for all $r>0$ there is a $(d, r, M_r)$-cover for $C_F(\Gamma)$.  As $\asdim C_F(\Gamma)$ is independent of $F$, we may simply write $\asdim\Gamma$.  When referring to the action of $\Gamma$ on itself, we mean the action by left multiplication.  }

\lemma{Suppose $\Gamma\acts X$ is a free action by a finitely generated group.  Then $\dad(\Gamma\acts X)\geq \asdim \Gamma$.  In fact, this inequality holds even if one modifies the definition of $\dad$ to allow covers by Borel sets, or even arbitrary sets.  }\label{easy inequality} 
\begin{proof}
Let $F$ be a finite generating set for $\Gamma$ as usual.  Suppose $\dad(\Gamma\acts X)\leq d$.  Let $r>0$.  Then there is a $(d, F^r, M_r)$-cover $\mathcal{V}$ for $\Gamma\acts X$.  Let $\mathcal{O}(x)$ be the orbit of some $x\in X$.  Identify $\mathcal{O}(x)$ with $C_F(\Gamma)$ by $\gamma\cdot x\mapsto \gamma$.  Then $\mathcal{V}$ gives a cover $\mathcal{V}'$ of $C_F(\Gamma)$.  Moreover, two points $y, y'\in \mathcal{O}(x)$ are connected by an $F^r$ chain iff their corresponding elements in $C_F(\Gamma)$ are connected by an $r$-chain.  Therefore, the $r$-components of each $V\in \mathcal{V}'$ are $M_r$-bounded, so $\mathcal{V}'$ is a $(d, r, M_r)$-cover for $C_F(\Gamma)$.  \end{proof}

\normalfont
Saying anything about the reverse inequality will require considerably more work.  


%
%

\normalfont


%
%

\normalfont

\section{Union theorem}

Union theorems are common in dimension theories.  The intuition is that a finite union of objects with dimension $\leq d$ should still have dimension $\leq d$.  Our method of showing this for $\dad$ comes from \cite[Lemma 3.5 and Corollary 3.6]{Brodskiy2006AHT}.  This section also includes an application of the union theorem which will be used later, as well as some other technical lemmas.  

%
\lemma{Suppose $\Gamma\acts X$ is free and $A, B\subset X$ such that the $F^r$-components of $B$ have cardinality at most $R$ and the $F^{r(R + 2)}$-components of $A$ have cardinality at most $D$.  Then the $F^r$-components of $A\cup B$ have cardinality at most $\#B_e^{(D-1)(R+1) + 1}(C_{F^r}(\Gamma))$.  }

\begin{proof}
Consider an $F^{r}$-chain in $A\cup B$: $x_0, \ldots, x_n$ and assume this chain has no repeated points.  Suppose $x_i$ and $x_j$ are two consecutive elements of $A$, so $x_{i+1}, \ldots, x_{j-1}$ are all points in $B$.  Then these points form an $F^r$-chain in $B$ which therefore has cardinality at most $R$, and therefore length at most $R$ (since there are no repeated points).  But then $x_{i+1}$ and $x_{j-1}$ are connected by something in $F^{rR}$, and so $x_i$ and $x_j$ are in the same $F^{r(R + 2)}$-component of $A$.  By similar reasoning, all elements of $A$ in the original chain are in the same $F^{r(R + 2)}$-component of $A$.  The subchain $x_0, \ldots, x_n$ of points contained in $A$ then has cardinality at most $D$, hence length at most $D$.  Since each chain between each element of $A$ has length at most $R$, the length of the original chain is at most $(D-1)(R+1) + 1$.  The $F^r$-components of $A\cup B$ therefore have cardinality at most $\#B_e^{(D-1)(R+1) + 1}(C_{F^r}(\Gamma))$.  \end{proof}


\theorem{(Finite Union Theorem) Suppose $\Gamma\acts X$ and $\{A_i\}_{i=0}^K$ are open subsets of $X$.  Suppose $f_i$ is a $d$-dimensional control function for $\Gamma\acts A_i$.  Let $r>0$.  Define $r_i$ and $R_i$ inductively starting with $r_0 = r$, $R_0 = f_0(r)$ and then for $1\leq i\leq K$ defining $r_{i} = r(R_{i-1} + 2)$, and $$R_{i} = \#B_e^{(R_{i-1} - 1)(f_i(r_{i}) + 1) + 1}(C_{F^r}(\Gamma)).$$  Then there is a $(d, F^r, R_K)$-cover for $\Gamma\acts \cup_i A_i$.  Notice that $R_K$ depends only on the $f_i$ and $K$.}  \label{union theorem}


\begin{proof}
Take a $(d, F^{r_i}, f_i(r_i))$-cover $\mathcal{V}^{(i)} =\{V_j^{(i)}\}_{j=0}^d$ of $A_i$ for each $i$.  Form a cover $\mathcal{V}$ of $\cup_{i=0}^K A_i$ by putting $\mathcal{V}_j = \cup_{i=0}^K V_j^{(i)}$.  Then use the previous lemma and induction.   At the $s$-th step of the induction (a union of $s+1$ sets), we have $r = r$ and $rR + 2r = rR_{s-1} + 2r = r_s$ ($s = 1, \ldots, K$).  \end{proof}


\normalfont

\definition{Suppose $\Gamma$ has finite generating set $F$ and $(X, \Gamma, \{D_\gamma\}_{\gamma\in \Gamma}, \{\theta_\gamma\}_{\gamma\in \Gamma})$ is a partial dynamical system.  For $P\in \N$, we define $$B_x^P(X, \Gamma, \{D_\gamma\}_{\gamma\in \Gamma}, \{\theta_\gamma\}_{\gamma\in \Gamma}) := \{y\in X \mid y = \theta_\gamma(x) \text{ for some } \gamma\in B_e^P(C_F(\Gamma))\}.$$  We may also just write $B_x^P$ or $B_x^P(X)$ when unambiguous.  }

\lemma{Suppose $(X, \Gamma, \{D_\gamma\}_{\gamma\in \Gamma}, \{\theta_\gamma\}_{\gamma\in \Gamma})$ is a partial dynamical system such that $\theta_\gamma(x) = \theta_{\gamma'}(x)$ for $\gamma, \gamma'\in B^P_e(C_F(\Gamma))$ if and only if $\gamma = \gamma'$.  Then this system restricted to any $B_x^P$ is conjugate to a restricted action $\Gamma\acts A\subset B_e^P(C_F(\Gamma))$ by left multiplication.  }\label{locally free lemma}

\begin{proof}
Consider the subset $B_x^P = \{y_\gamma \mid y_\gamma = \theta_\gamma(x) \text{ with } \gamma\in B_e^P(C_F(\Gamma))\}$.  Identify $B_x^P$ with a subset $A$ of $B_e^P(C_F(\Gamma))$ by the map $f(y_\gamma):= \gamma$.  Then if $y_\delta, y_{\delta'}\in B_x^P$ are such that $\theta_\gamma(y_\delta) = y_{\delta'}$, we have $y_{\delta'} = \theta_\gamma(y_\delta) = \theta_\gamma\circ \theta_\delta(x) = \theta_{\gamma\delta}(x) = y_{\gamma\delta}$ implying $\gamma\delta = \delta'$.  Hence, $f(\theta_\gamma(y_\delta)) = f(y_{\delta'}) = \delta' = \gamma\delta = \gamma\cdot \delta$, so $f$ is equivariant.  \end{proof} 

\lemma{Suppose $\Gamma$ is a finitely generated group with generating set $F$ and $C_F(\Gamma)$ is a Cayley graph.  Let $(X, \Gamma, \{D_\gamma\}_{\gamma\in \Gamma}, \{\theta_\gamma\}_{\gamma\in \Gamma})$ be a free partial dynamical system.  Then a $(d, r, M)$-cover for $C_F(\Gamma)$ gives rise to a $(d, F^r, \#B_e^M(C_F(\Gamma))$-cover for $(X, \Gamma, \{D_\gamma\}_{\gamma\in \Gamma}, \{\theta_\gamma\}_{\gamma\in \Gamma})$.  } \label{geometry to dynamics lemma}
\begin{proof} For each orbit of $(X, \Gamma, \{D_\gamma\}_{\gamma\in \Gamma}, \{\theta_\gamma\}_{\gamma\in \Gamma})$, identify the points in $X$ with a subset of $C_F(\Gamma)$ by choosing any $x_0$ to correspond to the identity and then identifying $\theta_\gamma(x)$ with $\gamma$.  Then notice that two points are connected by an $F^r$ chain iff their corresponding elements in $C_F(\Gamma)$ are connected by an $r$-chain.  \end{proof}

\normalfont

Next, we want to establish some lemmas about the $\dad$ of sequences of actions.  The content of the following lemma is essentially that of \cite[3.1]{boxspacesDT} and mirrors the theory of box spaces of residually finite groups.  In the proof of the main theorem, this is part of what allows us to infer a tight bound on the dimension whenever the dimension is finite.  

\lemma{Suppose $\Gamma$ is finitely generated with generating set $F$ and $(G_n, \Gamma, \{D^n_\gamma\}_{\gamma\in \Gamma}, \{\theta^n_\gamma\}_{\gamma\in \Gamma})_n$ is a sequence of partial dynamical systems with $G_n$ finite such that for all $r>0$ there is $N>0$ such that there is a $(d, F^r, M_r)$-cover for $(G_n, \Gamma, \{D^n_\gamma\}_{\gamma\in \Gamma}, \{\theta^n_\gamma\}_{\gamma\in \Gamma})$ for all $n\geq N$ ($M_r$ does not depend on $n$).  Suppose also that for all $P>0$ there is $N'$ such that for $n\geq N'$ and any $r'>0$, the system $(G_n, \Gamma, \{D^n_\gamma\}_{\gamma\in \Gamma}, \{\theta^n_\gamma\}_{\gamma\in \Gamma})_n$ restricted to any $B_x^P$ ($x\in G_n$) has a $(d', F^{r'}, M_{r'}')$-cover with $M'$ depending only on $r'$ and $\Gamma$.  Then for all $r>0$, there is $N''$ and a $(d', F^r, M_r'')$-cover for $(G_n, \Gamma, \{D^n_\gamma\}_{\gamma\in \Gamma}, \{\theta^n_\gamma\}_{\gamma\in \Gamma})_n$ for all $n\geq N''$.  Moreover, $M_r''$ depends only on $d$ and the function $r\mapsto M'_r$.  } \label{DT lemma generalized}

\begin{proof}Fix $r>0$.  For each $0\leq i\leq d$, define $f_i(r) = M_r'$ for all $i$.  Define $r_i$ and $R_i$ inductively as in \ref{union theorem}.  By hypothesis, there is $N$ such that for $n\geq N$, there exists a $(d, F^{R_K}, M_{R_K})$-cover $\mathcal{V} = \{V_i\}_{i=0}^{d}$ for $(G_n, \Gamma, \{D^n_\gamma\}_{\gamma\in \Gamma}, \{\theta^n_\gamma\}_{\gamma\in \Gamma})$.  

Then find $N'\geq N$ so that $(G_n, \Gamma, \{D^n_\gamma\}_{\gamma\in \Gamma}, \{\theta^n_\gamma\}_{\gamma\in \Gamma})_n$ restricted to any $B_x^{M_{R_K} + 1}$ has a $(d', F^{r_i}, f_i(r_i))$-cover.  Find such a cover $\mathcal{U}^i = \{U_k^i\}_{k=0}^{d'}$ for each $F^{R_K}$-component of $V_i$.  For each $0\leq i\leq d$, form a $(d', F^{r_i}, f_i(r_i))$-cover $\mathcal{U} = \{U_k\}_{k=0}^{d'}$ of $V_i$ by putting $U_k = \bigcup_j U_k^j$.  It follows from \ref{union theorem} that there is a $(d', F^{r}, g_d(r))$-cover for each $(G_n, \Gamma, \{D^n_\gamma\}_{\gamma\in \Gamma}, \{\theta^n_\gamma\}_{\gamma\in \Gamma})$ with $n\geq N'$.  \end{proof}

\definition{If $\Gamma$ is a group and $C_F(\Gamma)$  is a Cayley graph ($F$ is a finite generating set), we say $\Gamma$ has polynomial growth of degree at most $d$ if there are $d>0$ and $C>0$ such that if $\mathcal{G}(r) = \#B_e^r(C_F(\Gamma))$, then $\mathcal{G}_n(r)\leq Cr^d$.  A sequence of partial dynamical systems $(G_n, \Gamma, \{D^n_\gamma\}_{\gamma\in \Gamma}, \{\theta^n_\gamma\}_{\gamma\in \Gamma})_n$ has polynomial growth of degree at most $d$ if there are $d>0$ and $C>0$ such that if $\mathcal{G}_n(r) = \sup_{x\in G_n} \#\{y : y = \gamma\cdot x\ \text{ for } \gamma\in B^r_e(C_F(\Gamma))\}$, then $\mathcal{G}_n(r)\leq Cr^d$ for all $n$.  }\label{polynomial growth definition}

We will make use of a standard `greedy algorithm' in this section and again later, so we formalize such an algorithm with a lemma.  

\lemma{Suppose $S$ is a finite set and $\sim$ is a reflexive, symmetric relation on $S$ such that $\#\{t\in S \mid t\neq s\text{ and } s\sim t\}\leq D$ for all $s\in S$.  Then there is a cover $\mathcal{C} = \{C_i\}_{i=0}^D$ such that $s\nsim t$ for all $s\neq t$ in $\mathcal{U}_i$ for all $i$.  }
\begin{proof} Assign elements to the sets in $\mathcal{C}$ in the following way.  Start by picking any element $s_0\in S$.  Then there are at most $D$ other elements of $S$ which are related to $s_0$.  Assign $s_0\in C_0$ and those other elements to $C_1$ through $C_D$ so that at most one belongs to each.  Then consider one of the elements just discussed other than $s_0$.  This element is related to at most $D$ other elements of $S$ including $s_0$, so we can repeat the process.  Continue in this way until all elements of $S$ are assigned to some set in $\mathcal{C}$.  \end{proof} \label{greedy algorithm}

\lemma{Suppose $\Gamma$ is finitely generated and $(G_n, \Gamma, \{D^n_\gamma\}_{\gamma\in \Gamma}, \{\theta^n_\gamma\}_{\gamma\in \Gamma})_n$ is a sequence of partial dynamical systems on finite sets with polynomial growth of degree at most $d$.  Then $\dad(G_n, \Gamma, \{D^n_\gamma\}_{\gamma\in \Gamma}, \{\theta^n_\gamma\}_{\gamma\in \Gamma})_n\leq K = 4^d + 1$.  }\label{polynomial growth lemma}
\begin{proof}
The proof is essentially the same as that of \cite[3.3]{boxspacesDT}.  

Fix $F\subset \Gamma$ a finite generating set.  For a partial dynamical system $(G, \Gamma, \{D_\gamma\}_{\gamma\in \Gamma}, \{\theta_\gamma\}_{\gamma\in \Gamma})$ on a discrete space, the following condition implies its dynamic asymptotic dimension is at most $K-1$: for all $R>0$ there is $M>0$ and a cover $\mathcal{V}$ of $G$ by sets of cardinality at most some $M$ such that if $x\in V\in\mathcal{V}$, there are at most $K-1$ other elements $V'\in\mathcal{V}$ such that $\gamma\cdot x\in V'$ for some $\gamma\in F^R$.  To see this, notice that one can use a greedy algorithm as in \ref{greedy algorithm} to sort the elements of $\mathcal{V}$ into $K$ subcollections such that taking the union over each subcollection yields a $(d, F^r, M)$-cover.  

By assumption, we have that $\#B_x^R(G_n)\leq CR^d$ for all $R>0$, $x\in G_n$ for any $n$.  Let $K = 4^d + 1$, $R> 1$ and $S_0 = 4^{m+1}R$ where $m$ is such that $(K/4^d)^m\geq CR^d$.  We claim that for every $n$ there exists $R_n$ such that $R\leq R_n\leq \frac{S_0}{4}$ and $|B^{4R_n}_x(G_n)|\leq K|B^{R_n}_x(G_n)|$ (where $x\in G_n$ is arbitrary).  If this were not the case then $K^i|B_x^R(G_n)|<|B^{4^iR}_x(G_n)|\leq C4^{id}R^d$ for every $i$ with $4^iR\leq S_0$, so setting $i = m$ would give $C4^{md}R^d>K^m|B^R_x(G_n)|\geq C4^{md}R^d|B^R_x(G_n)|$.  But then we would have $1>|B^R_x(G_n)|$, a contradiction.  

Now, take $X_n$ maximal in $G_n$ such that $B_x^{R_n}(G_n)$ and $B_y^{R_n}(G_n)$ are disjoint for all $x, y\in X_n$.  Consider the collection $\mathcal{V} = \{B_x^{2R_n}(G_n) \mid x\in X_n\}$.  If $z\in G_n$, maximality of $X_n$ implies there is $x\in X_m$ such that $B_z^{R_n}(G_n)\cap B_x^{R_n}(G_n)$ and so $z\in B_x^{2R_n}(G_n)$, and so $\mathcal{V}$ covers $G_n$.  

Finally, we check that if $x\in G_n$ and $x\in V\in\mathcal{V}$, there are at most $K$ other elements $V'\in\mathcal{V}$ such that $\gamma\cdot x\in V'$ for some $\gamma\in F^{R_n}$.  Let $z\in G_n$.  For every $B_x^{2R_n}(G_n)\in \mathcal{V}$ which has an element at a distance at most $R_n$ to $z$, we have that $x\in B_z^{2R_n + R_n}(G_n)\subset B_z^{3R_n}(G_n)$.  Now consider $B_z^{3R_n}(G_n)\cap X_n$.  Since $B_x^{R_n}(G_n)$ and $B_y^{R_n}(G_n)$ are disjoint for any $x$ and $y$ in $X_m$, we have that $\#(B_z^{3R_n}(G_n)\cap X_n)\leq \frac{\#B_z^{4R_n}(G_n)}{\#B_z^{R_n}(G_n)}\leq 4^d + 1$.  

If $V'\in \mathcal{V}$ is such that there is $\gamma\in F^{R_n}$ with $\gamma\cdot V \cap V' \neq \emptyset$, then there is $z\in V'$ with $z\in B_z^{3R_n}(G_n)\cap X_n$.  Therefore, there are at most $4^d + 1$ such $V'\in \mathcal{V}$.  This completes the proof.  \end{proof}

\lemma{Suppose $\Gamma$ has polynomial growth and $(X, \Gamma, \{D_\gamma\}_{\gamma\in \Gamma}, \{\theta_\gamma\}_{\gamma\in \Gamma})$ is a partial dynamical system.  Then $(X, \Gamma, \{D_\gamma\}_{\gamma\in \Gamma}, \{\theta_\gamma\}_{\gamma\in \Gamma})$ has polynomial growth bounded by that of $\Gamma$.  That is, a $d>0$ and $C>0$ that work for $\Gamma$ in the sense of \ref{polynomial growth definition} also work for the whole sequence of partial actions.  } \label{polynomial growth lemma}



%
%
%

\section{Residual finiteness}
\normalfont
We will also make use of a property which allows group actions to be approximated by actions on finite sets.  The definition below was originally introduced by Kerr and Nowak in \cite{kerr2011residually}.  

\definition{An action $\Gamma\acts X$ on a metric space by homeomorphisms is residually finite if for all $F\subset \Gamma$ finite and $\epsilon>0$ there is a finite set $E$, an action $\Gamma\acts E$, and a map $\zeta: E\to X$ such that $\zeta(E)$ is $\epsilon$-dense in $X$ (meaning every $\epsilon$-ball in $X$ intersects $E$ nontrivially) and that $d(\zeta(\gamma\cdot e), \gamma \cdot\zeta(e))<\epsilon$ for all $\gamma\in F$.  Such an action is called a $(F, \epsilon)$-approximating action for $\Gamma\acts X$.  If $X$ has no isolated points, a perturbation argument shows we can assume $\zeta$ to be an inclusion.  If $(X, \Gamma, \{D_\gamma\}_{\gamma\in \Gamma}, \{\theta_\gamma\}_{\gamma\in \Gamma})$ is a partial dynamical system, $E$ is a finite set and $\zeta: E\to X$ is a map, we similarly say a partial dynamical system $(E, \Gamma, \{D^E_\gamma\}_{\gamma\in \Gamma}, \{\theta^E_\gamma\}_{\gamma\in \Gamma})$ is $(F, \epsilon)$-approximating if $E$ is $\epsilon$-dense in $X$ and if, whenever $\theta_\gamma(\zeta(x))$ is defined for $\gamma\in F$ and $x\in E$, we have that $\theta^E_\gamma(x)$ is defined and $d(\theta_\gamma(\zeta(x)), \zeta(\theta_\gamma^E(x)))<\epsilon$.   }

\normalfont
As with the theory of residually finite groups, there is a close relationship between isometry and residual finiteness.  We can make use of residual finiteness of actions later in \ref{asymptotically faithful lemma} without loss of generality in the main theorem thanks to the following result.  

\theorem[\cite{pilgrim_2022}*{Theorem 3.6}]{Every faithful, isometric action by a finitely generated, amenable group is residually finite.  }\label{RF theorem} \qed

\normalfont
The work of \cite{boxspacesDT} shows that the asymptotic dimension of box spaces of residually finite groups is closely related to the asymptotic dimension of the group itself.  As the dynamic asymptotic dimension behaves somewhat like the asymptotic dimension of a box space (both are infinite in the non-amenable case for example), we expect that residual finiteness of actions may be helpful to relating their DAD to the asymptotic dimension of the acting group.  

\section{Cantor-like decompositions}

\normalfont
A key ingredient in the proof of the main theorem is a technical condition on $X$ which allows an isometric, residually finite action to be modeled by several sequences of partial dynamical systems on finite spaces, thereby allowing us to apply the work of the previous section.  We will see that many spaces satisfy this property.  

Isometric actions on Cantor sets can be described as inverse limits of finite actions \cite[2.3]{pilgrim_2022}, which allows their $\dad$ to be studied more easily.  We therefore would like to cover $X$ by Cantor sets and use the union theorem from section 3.  Instead of doing exactly this, we cover $X$ by sets which (as $\epsilon\to 0$) resemble the sets one might intersect in constructing a Cantor set-like object.  

\definition{Let $X$ be a compact metric space and let $\epsilon, \delta>0$.  We say $X$ is $(K, \epsilon, \delta)$-Cantor decomposable if there is an open cover $\mathcal{U} = \mathcal{U}_0 \sqcup \ldots \sqcup \mathcal{U}_K$ of $X$ with $\diam(U)<\epsilon$ for all $U\in \mathcal{U}_i$ and all $i$, and $d(U, V)>\delta$ for all distinct $U, V\in \mathcal{U}_i$ and all $i$.  We call such a cover a $(K, \epsilon, \delta)$-decomposition.  }

\normalfont

This condition comes from trying to force an argument like \cite[2.3]{pilgrim_2022} to work by dividing $X$ into different subspaces.  Notice that the definition is somewhat reminiscent of asymptotic dimension, but small-scale.  We therefore expect $K$ to be related to some notion of metric dimension for $X$.

\normalfont

\definition{A metric space $(X, d)$ is \textit{$M$-doubling} if for all $r>0$ and $x\in X$, the closed ball $B_x^{2r}$ can be covered by at most $M$ closed balls of radius $r$.  Further define $\text{dim}_{d}(X) := \log_2M$, which we call the \textit{doubling dimension} of $X$.  } 

\normalfont
The origins of this concept are not precisely clear and it may have been rediscovered independently multiple times.  It can be reasonably attributed to Assouad, though his is a different, equivalent definition.  See the introduction of \cite{doublingspaces} for additional exposition.  

\example{The space $\mathbb{R}$ is $2$-doubling and therefore has doubling dimension $1$.  Working out $M$ exactly for higher-dimension Euclidean spaces is difficult, but it is more straightforward to see that $\mathbb{R}^d$ is $M$-doubling for \textit{some} $M\leq 4\cdot2^d$.  This follows since a $d$-cube with side length $l$ can be covered by $4$ $d$-balls of diameter $l$, and it takes $2^d$ $d$-cubes with side length $l/2$ to cover one with side length $l$.  More generally, a Riemannian manifold admits an isometric embedding into some Euclidean space by the Nash embedding theorem, and therefore has finite doubling dimension.  } \label{doubling example}

\normalfont
As alluded to earlier, we will now see that Cantor decomposability relates quantitatively to the doubling dimension.  

\lemma{Suppose $(X, d)$ is compact and $M$-doubling.  Then $(X, d)$ has a $(M^{\lceil2 + \log_2(k+3)\rceil},\epsilon, k\epsilon)$-decomposition for all $\epsilon>0$ and $k\in \N$.  If $X$ has a $(K, \epsilon, 2\epsilon)$-decomposition for all $\epsilon>0$, then $X$ is $K$-doubling.  In particular, this implies $X$ has finite doubling dimension iff for all $k$ there exists $K$ such that $X$ is $(K, \epsilon, k\epsilon)$-Cantor decomposable for all $\epsilon>0$.  } \label{doubling dimension lemma}

\begin{proof}
For the first claim, fix $\epsilon>0$.  Let $\mathcal{C}$ be a cover of $X$ by (closed) $\epsilon/2$-balls using as few balls as possible.  If $x\in X$ and $\mathcal{U} = \{B\in \mathcal{C}| B\cap B_x^{(k+2)\epsilon} \neq \emptyset\}$ then $\mathcal{U}$ is a cover of $\bigcup_{B\in \mathcal{U}} B$ using as few $\epsilon/2$-balls as possible: if $\bigcup_{B\in \mathcal{U}} B$ had a cover $\mathcal{U}'$ by fewer balls, then $(\mathcal{C}\setminus \mathcal{U})\cup\mathcal{U}'$ would be a cover of $X$ by fewer balls.  Therefore, since $\bigcup_{B\in \mathcal{U}} B\subset B_x^{(k+3)\epsilon}$, $\#\mathcal{U}\leq M^{\lceil2 + \log_2(k+3)\rceil}$ (where $\lceil\cdot\rceil$ denotes the ceiling function).  

It is therefore possible to color the balls using $\#\mathcal{U}\leq M^{\lceil2 + \log_2(k+3)\rceil}$ colors such that any two balls of the same color are distance at least $(k+1)\epsilon$ apart.  This can be shown via a standard `greedy algorithm' (see \ref{greedy algorithm}).  Finally, we can replace each ball with an open neighborhood of itself while keeping the diameter of each $<\epsilon$ and keeping balls in the same collection $>k\epsilon$ apart.  

For the second claim, suppose $X$ is $(K, \epsilon, 2\epsilon)$-Cantor decomposable for all $\epsilon>0$.  Given a ball $B\subset X$ of radius $2\epsilon$, find a cover of $X$ witnessing its $(K, \epsilon, 2\epsilon)$-Cantor decomposability.  Then the subcollection of sets in this cover which intersect $B$ is a cover of $B$.  Moreover, any set intersecting $B$ is within $2\epsilon$ of any other such set, so this cover of $B$ contains at most $K$ sets of diameter $<\epsilon$.  \end{proof}

\normalfont
The existence of $(K, \epsilon, \delta)$-decompositions for $\delta$ large compared to $\epsilon$ really does require finite doubling dimension rather than finiteness of other types of dimension (e.g. Lebesgue covering dimension).  We give an example to show this.  

\example{There exists a Cantor set with infinite doubling dimension.  }

\begin{proof}
Let $\mathcal{C}$ be the space of sequences $(a_n)_{n=1}^\infty$ where $a_n\in \{1, \ldots, n\}$.  Give $\mathcal{C}$ the metric $d((a_n), (b_n)) = \sum_{n} \frac{|a_n - b_n|}{n2^n}$.  Then a ball of radius $\frac{1}{2^k}$ centered at $(a_n)$ consists of sequences which agree with $(a_n)$ up to at least index $k$.  Such a ball therefore requires $k+1$ balls of radius $\frac{1}{2^{k+1}}$ to cover it.   \end{proof}

\normalfont
We conclude this section by showing how Cantor decompositions can be used to model certain isometric actions by partial dynamical systems acting on decompositions of $X$.  

\lemma{Suppose $\Gamma\acts X$ (denoted by $\gamma\cdot x$) is isometric and residually finite, that $\Gamma$ has finite generating set $F$, and that $X$ has no isolated points and finite doubling dimension.  Let $r>0$ and $P>0$.  Then there is $\epsilon_0$ such that if $\epsilon<\epsilon_0$ and $\mathcal{U}^\epsilon = \mathcal{U}_0^\epsilon \sqcup \cdots \sqcup \mathcal{U}_K^\epsilon$ is a $(K, \epsilon, 15\epsilon)$-decomposition of $X$ such that, for all $i$ and $U\in \mathcal{U}_i$, $U$ contains a ball of radius $\epsilon/2$ (this can be done by first taking a $(K, \epsilon/2, 16\epsilon)$-decomposition for each $\epsilon$ and then replacing each set by its $\epsilon/2$-neighborhood); then there is a partial dynamical system $(\mathcal{U}_i^{\epsilon}, \Gamma, \{D_\gamma\}_{\gamma\in \Gamma}, \{\theta_\gamma\}_{\gamma\in \Gamma})$ such that

\begin{description}
\item(i) If $\gamma\in F^r$ and $\gamma\cdot U\cap V\neq \emptyset$ for $U, V\in \mathcal{U}_i^{\epsilon}$, then $\theta_\gamma(U) = V$.  


\item(ii) If $w = f_j\cdots f_1$ and $w' = f'_p\cdots f'_1$ are words in $F^r$ with length $\leq P + 1$ which are not equal as elements of $\Gamma$, then $\theta_{f_j}\circ \cdots\circ \theta_{f_1}(x) \neq \theta_{f'_p}\circ\cdots\circ \theta_{f'_1}(x)$ for all $x$ where both sides are defined.

\end{description}

}\label{asymptotically faithful lemma}

\begin{proof}
Start by letting $\epsilon$ be sufficiently small that translates (according to the action $\Gamma\acts X$) of subsets of diameter $<\epsilon$ by elements of $F^{(P+1)r}$ are disjoint (using that $\Gamma\acts X$ is free).  This will later be used to ensure condition (ii) holds.  

Let $\Gamma\acts E\subset X$ be a $(F^r, \epsilon/2)$-approximating action for $\Gamma\acts X$.  Use $\delta_\gamma(x)$ to denote the action $\Gamma\acts E$.  

Now, for each $i$, we can restrict the action $\Gamma\acts E$ to a partial action $\Gamma\acts E\cap N_{2\epsilon}(\bigcup_{U\in \mathcal{U}_i^\epsilon}U)$ (where $N_{2\epsilon}(\cdot)$ denotes the open $2\epsilon$-neighborhood).  Each $x\in E\cap N_{2\epsilon}(\bigcup_{U\in \mathcal{U}_i^\epsilon} U)$ can be associated to an $\tilde{x}\in \bigcup_{U\in \mathcal{U}_i^\epsilon} U$ in such a way that $\tilde{x} = x$ if $x\in \bigcup_{U\in \mathcal{U}_i^\epsilon} U$, $d(x, \tilde{x})<3\epsilon$, and $\tilde{x}\neq\tilde{y}$ for $x\neq y$ (since $X$ has no isolated points).  Let $\tilde{E} = \{\tilde{x}\ \mid x\in E\cap N_{2\epsilon}(\bigcup_{U\in \mathcal{U}_i^\epsilon} U)\}$ and conjugate by the correspondence $x\mapsto \tilde{x}$ to obtain a partial dynamical system $(\tilde{E}, \Gamma, \{\alpha_\gamma\}_{\gamma\in \Gamma}, \{C_\gamma\}_{\gamma\in \Gamma})$.  

Then, if there is $\gamma\in F^r$ and $\gamma\cdot U\cap V\neq \emptyset$ for some $V\in \mathcal{U}_i^{\epsilon}$; then (as $E$ is $\epsilon/2$-dense and as $U$ contains an $\epsilon/2$-ball) there is $x\in U$ with $x\in E$ so that $\tilde{x} = x$ and $d(\delta_\gamma(x), V)<2\epsilon$, and so $\widetilde{\delta_\gamma(x)} = \alpha_\gamma(x)\in V$.  This will guarantee condition (i).  Also notice that $d(\gamma\cdot \tilde{x}, \alpha_\gamma(\tilde{x}))<7\epsilon$ for any $\tilde{x}\in \tilde{E}$.  

Then the partial dynamical system $(\tilde{E}, \Gamma, \{\alpha_\gamma\}_{\gamma\in \Gamma}, \{C_\gamma\}_{\gamma\in \Gamma})$ is $(F^r, 7\epsilon)$-approximating for the restricted partial action $\Gamma\acts \bigcup_{U\in \mathcal{U}_i^\epsilon}U$ and every $U\in \mathcal{U}_i^\epsilon$ contains at least one element of $\tilde{E}$ (since the set $E$ was $\epsilon/2$-dense).  

Moreover, if $x\in U$ with $x\in \tilde{E}$ and $\alpha_\gamma(x)\in V\in \mathcal{U}_i^\epsilon$, and $y\in U$ with $y\in \tilde{E}$; then if $\alpha_\gamma(y)\in \bigcup_{U\in \mathcal{U}_i^\epsilon}U$, $d(\alpha_\gamma(x), \alpha_\gamma(y))<15\epsilon$, and so $\alpha_\gamma(y)\in V$.  For $\gamma\in \Gamma$, restrict $\alpha_\gamma$ to the points $x\in \tilde{E}$ such that there exists a word $f_m\cdots f_1$ in $F^r$ which is equal to $\gamma$ in $\Gamma$ and such that the composition $\alpha_{f_m}\circ \cdots\circ\alpha_{f_1}(x)$ is defined (recall that since we are restricting a partial action of $\Gamma$, such compositions must agree with $\alpha_\gamma$ whenever they are defined).  

Having been restricted this way, each partial bijection $\alpha_\gamma$ for $\gamma\in \Gamma$ descends to a partial bijection $\beta_\gamma$ of $\mathcal{U}_i^\epsilon$, and form the partial dynamical system $(\mathcal{U}_i^{\epsilon}, \Gamma, \{D_\gamma\}_{\gamma\in \Gamma}, \{\beta_\gamma\}_{\gamma\in \Gamma})$.  We have that $\beta_\gamma\circ \beta_\lambda\subseteq \beta_{\gamma\lambda}$ since $\beta_\gamma$ is induced from $\alpha_\gamma$ on a quotient space.  Explicitly, this definition amounts to saying $\beta_\gamma(U) = V$ iff there exists $\tilde{x}\in \tilde{E}\cap U$ and a sequence $f_1, \ldots, f_n\in F^r$ such that $\alpha_{f_m}\circ\cdots\circ \alpha_{f_1}(x)$ is defined and in $V$.  This does not depend on $\tilde{x}$ by the preceding paragraph, and does not depend on $f_1, \ldots, f_m$ since $\alpha_{f_m}\circ \cdots \circ \alpha_{f_1} = \alpha_\gamma$ whenever the left hand side is defined.  The descended self maps are still partial bijections since non-injectivity would imply $d(\gamma\cdot x, \gamma\cdot y)<15\epsilon$ for $x\in U\in \mathcal{U}_i^\epsilon$, $y\in V\in \mathcal{U}_i^\epsilon$, and $U\neq V$; and this would contradict that the original action is isometric.  As noted in the third paragraph, these systems satisfy (i).  

Now we make a further restriction.  For each $\gamma\in \Gamma$, replace $D_\gamma$ with the subset of itself consisting of those $U\in \mathcal{U}_i^{\epsilon}$ such that there exists a word $f_m\cdots f_1$ in $F^r$ which is equal to $\gamma$ in $\Gamma$ and such that for all $m_0\in \{1, \ldots, m\}$ there exists $V\in \mathcal{U}_i^\epsilon$ such that $\gamma_{m_0}\cdots \gamma_1\cdot U\cap V\neq \emptyset$ and let $\theta_\gamma$ be the restriction of $\beta_\gamma$ to this new $D_\gamma$.  Notice that this preserves (i) while now ensuring (ii) by the choice of $\epsilon$ we made at the beginning.  \end{proof}

\normalfont

Condition (ii) is needed since the partial dynamical systems constructed here are not free, so it is important that they be asymptotically free in an appropriate sense.  This condition is also used with \ref{locally free lemma} to show these partial dynamical systems locally resemble $\Gamma\acts \Gamma$.  Condition (i) ensures a $(F^{r}, M)$-cover for the system on $\mathcal{U}_i^\epsilon$ gives rise to a $(F^r, M)$-cover for the restricted action $\Gamma\acts \bigcup_{U\in \mathcal{U}_i^\epsilon} U$.  

The sequences of actions constructed above are therefore similar to box spaces of residually finite groups.  Indeed, we are motivated by past results relating dynamic asymptotic dimension to the asymptotic dimension of box spaces \cite{szabo2017rokhlin}.

\section{Application to DAD}

\lemma{Let $\Gamma\acts X$ and suppose there is a $(d, F^k, D)$-cover for $\Gamma\acts X$.  Then there is $\delta>0$ such that for $\mathcal{U}^\epsilon$ as in \ref{asymptotically faithful lemma} with $\epsilon<\delta$, $r = k$, and $p = D + 1$, there is a $(d, F^r, D)$-cover for $(\mathcal{U}_i^{\epsilon}, \Gamma, \{D_\gamma^{i, \epsilon}\}_{\gamma\in \Gamma}, \{\theta_\gamma^{i, \epsilon}\}_{\gamma\in \Gamma})$ for each $i$. }\label{bounded below lemma}

\begin{proof}
Fix $r>0$ and let $\mathcal{V} = \{V_m\}_{m=0}^d$ be a $(d, F^r, D)$-cover of $X$.  This cover has some Lebesgue number $\lambda$.  For $\epsilon<\lambda/4$ each element of $\mathcal{U}^{\epsilon}_i$ is entirely contained in at least one $V_m$, so we can form a cover of $\mathcal{U}_i$ by $d+1$ sets (by making some choices).  Moreover, for each $x\in U\in \mathcal{U}_i$ with $U\subset V_m$ for some $m$, any $F^r$-chain $x = x_0, \ldots, x_n$ with cardinality more than $D$ must have some $x_j\notin V_m$.  Replace each $V_m$ by its $\lambda/3$-interior, which can be done for all $i$ while maintaining that the $V_m$ cover $X$ and that each $U\in \mathcal{U}_i$ is entirely contained in some $V_m)$.  Now any $F^r$-chain as before with cardinality more than $D$ must have some $x_j\notin N_{\lambda/3}(V_m)$.  Then by taking $\epsilon<\frac{\lambda}{3(D+1)}$, any $F^r$-chain in $\mathcal{U}_i$ with more than $D$ elements beginning at $U\subset V_m$ will contain some $U'$ entirely outside of $V_m$.  \end{proof}

\theorem{Suppose $\Gamma$ is finitely generated, that $\Gamma\acts X$ is free and isometric, and that $X$ has finite doubling dimension and no isolated points.  Suppose further that $\dad(\Gamma\acts X)\leq d$ for some $d$.  Then  $\dad(\Gamma\acts X)\leq \asdim \Gamma$.  More precisely, if $F\subset \Gamma$, then for every $r>0$ there is a $(\asdim\Gamma, F^r, D)$-cover for $\Gamma\acts X$ and $D$ depends only on $K, r, C_F(\Gamma)$, and $d$.  } \label{og main theorem}

\begin{proof}
As always, assume $\Gamma$ has finite generating set $F$.  If $\Gamma$ is not amenable, then $\dad(\Gamma\acts X) = \infty$.  This follows from \cite[Corollary 8.27]{dasdimGWY}, which shows that a groupoid with finite dynamic asymptotic dimension must be amenable.  Since amenable measure-preserving actions must be by amenable groups (and isometric actions fix measures), this means the dimension is infinite when $\Gamma$ is non-amenable.  If $\Gamma$ is amenable, $\Gamma\acts X$ is residually finite by \ref{RF theorem}, and we can therefore apply \ref{asymptotically faithful lemma}.  

For each $0\leq i\leq K$, let $n = 1, 2, \ldots$, and form a sequence $(\mathcal{U}_i^{\epsilon_n}, \Gamma, \{D^{i, \epsilon_n}_\gamma\}_{\gamma\in \Gamma}, \{\theta^{i, \epsilon_n}_\gamma\}_{\gamma\in \Gamma})_n$ of partial dynamical system with the properties from \ref{asymptotically faithful lemma} with $r = n$, $P = n$ and $\epsilon_n < 1/n$.  

Since $\Gamma\acts X$ is finite dimensional, \ref{bounded below lemma} shows we have, for any $r>0$, a $(d, F^r, M_r)$-cover for these systems for $n$ sufficiently large and all $i$.  Then by \ref{asymptotically faithful lemma}(ii), \ref{locally free lemma}, \ref{geometry to dynamics lemma}. and \ref{DT lemma generalized}, we have a $(\asdim\Gamma, F^r, M'_r)$-cover for each of these systems and $M'_r$ depends only on $r$, $C_F(\Gamma)$, and $d$.  

Then for $n\geq r$, by \ref{asymptotically faithful lemma}(i), this gives a $(\asdim\Gamma, F^r, M'_r)$-cover for the restricted action $\Gamma\acts \bigcup_{U\in \mathcal{U}_i} U$ for each $i$.  


Define $f_i(r):=M'_r$ for all $i$ and let $r_i$ and $R_i$ be as in \ref{union theorem}.  Then, letting $n\geq R_K$, we can apply \ref{union theorem} to get a $(\asdim\Gamma, F^r, R_K)$-cover for $\Gamma\acts X$, while also ensuring $n$ is sufficiently large that \ref{asymptotically faithful lemma}(ii) holds with $P = R_K$.  Since $R_K$ depends only on $K$, $r$, $C_{F^r}(\Gamma)$, and $d$, we are done.  \end{proof}





\corollary{Suppose $\Gamma\acts X$ is free and isometric and that $X$ has finite doubling dimension and no isolated points (e.g. is a Riemannian manifld).  Then either $\dad(\Gamma\acts X) = \infty$ or $\dad(\Gamma\acts X) = \asdim\Gamma$.  }\label{main theorem}

\begin{proof}
This follows from \ref{og main theorem} and \ref{easy inequality}.  \end{proof}

\section{Further discussion}

\normalfont
It is now possible to give a complete description of the asymptotic dimension of translation actions on compact Lie groups in terms of the amenability and asymptotic dimension of the acting group.  

\theorem{Let $\Gamma\acts G$ be a translation action by a finitely generated subgroup of a compact Lie group.  Then $\dad(\Gamma\acts G) = \asdim\Gamma$ if $\Gamma$ is amenable and $\dad(\Gamma\acts G) = \infty$ otherwise.  }

\begin{proof} If $\Gamma$ is not amenable, then $\dad(\Gamma\acts G) = \infty$ by the same argument as in \ref{og main theorem}, as this action is Haar measure-preserving.  Assume therefore that $\Gamma$ is amenable.  

Representation theory shows $G\leq U_n$ is a subgroup of the $n\times n$ unitaries for some $n$.  Since $\dad$ does not increase when passing to a closed, invariant subset, we can assume $G = U_n$. 

Fix a finite subset $F\subset \Gamma$.  The group $\Gamma$ is a finitely generated, amenable subgroup of a compact Lie group, and is therefore virtually abelian.  This follows from a combination of the Tits alternative, Lie's theorem, and Engel's theorem (see, for instance, \cite{stack} for a proof).  

In particular, this means $\Gamma$ is virtually nilpotent, hence polynomial growth.  Moreover, $U_n\subset M_n(\C)\cong \C^{n^2}$ has finite doubling dimension.  

The sequence of partial dynamical systems constructed in \ref{asymptotically faithful lemma} have uniform polynomial growth by \ref{polynomial growth lemma} and so have asymptotic dimension uniformly $\leq d$ for some $d$.  We can therefore apply \ref{main theorem}.  \end{proof}

\normalfont
We can similarly compute the $\dad$ of many actions by $\Z^d $.  

\theorem{Suppose $\Z^d\acts X$ is a free, isometric action and $X$ has finite doubling dimension (e.g. is a Riemannian manifold).  Then $\dad(\Z^d\acts X) = d$.  }

\begin{proof}
Apply \ref{og main theorem} and \ref{polynomial growth lemma} as in the proof above.  \end{proof}

\normalfont
Although previously suspected, it remains somewhat surprising that the dynamic asymptotic dimension does not depend on the topological complexity of $X$.  Said another way, \ref{main theorem} shows that dynamic asymptotic dimension can be equivalently defined using covers by Borel sets or even arbitrary sets (at least for actions satisfying the hypotheses of \ref{main theorem}).  This is spiritually similar to the results about Borel asymptotic dimension found in \cite[Theorems 4.2 and 6.4]{conley2020borel}.    

We can also obtain more precise bounds for other dynamical dimensions.  As mentioned earlier, the dynamic asymptotic dimension is related to other dimension theories for group actions by \cite[5.1.4]{kerr2017dimension}.  This and corollary \ref{main theorem} show together that if $\dim(\Gamma\acts X)$ is the amenability dimension (or the tower dimension), then either $\dim(\Gamma\acts X)=\infty$ or $\dim(\Gamma\acts X)\leq (\asdim\Gamma + 1)(\dim X + 1) - 1$ where $\dim X$ is the covering dimension of $X$.

\bibliography{mybibliography2.bib}
\bibliographystyle{plain}

\end{document}